\documentclass[12pt]{article}
\usepackage{amsfonts}
\usepackage{amssymb}
\usepackage{amsmath}
\usepackage{graphicx}
\providecommand{\U}[1]{\protect\rule{.1in}{.1in}}
\newtheorem{theorem}{Theorem}

\newtheorem{remark}[theorem]{Remark}
\newcommand{\R}{\mathbb{R}}
\newcommand\Sp{\mathbb{S}}
\newcommand\Hy{\mathbb{H}}
\newcommand\Q{\mathbb{Q}}
\newcommand\bea{\begin{eqnarray*}}
\newcommand\eea{\end{eqnarray*}}
\newcommand\be{\begin{equation}}
\newcommand\ee{\end{equation}}
\newcommand{\po}{{\hspace*{-1ex}}{\bf .  }}

\def\<{\langle}
\def\>{\rangle}
\newcommand\qed{\ifhmode\unskip\nobreak\fi\ifmmode\ifinner\else
\hskip5 pt \fi\fi\hbox{\hskip5 pt \vrule width4 pt height6 pt
depth1.5 pt \hskip 1pt }}
\begin{document}

\title{An interior gradient estimate for the mean curvature equation of Killing
graphs and applications}
\author{M. Dajczer\thanks{Partially supported by CNPq and
FAPERJ.},\,\,J. H. de Lira\thanks{Partially supported by CNPq and
FUNCAP/PRONEX.}\,\,\,\,and\, J. Ripoll\thanks{Partially supported by CNPq.}}
\date{}
\maketitle

\begin{abstract}
We extend the interior gradient estimate due to Korevaar-Simon \cite{K} for solutions
of the mean curvature equation from the case of Euclidean graphs 
to the general case of Killing graphs. Our main application is the proof of existence
of Killing graphs with prescribed mean curvature function for continuous
boundary data, thus extending a result due to  Dajczer, Hinojosa and Lira 
\cite{DHL}. In addition, we prove the existence and uniqueness of radial
graphs in hyperbolic space with prescribed mean curvature function and
asymptotic boundary data at infinity.
\end{abstract}

\section{Introduction}

Gradient  bounds are fundamental a priori estimates for elliptic and
parabolic equations and play a key role in geometry and PDE. One of the most
important application of interior gradient estimates is the proof due to Serrin \cite{Se}
of the solvability of the Dirichlet problem for the constant mean curvature
hypersurface equation in a bounded $C^2$ domain of the Euclidean space for
just continuous boundary data; see Theorem 16.11 of \cite{GT}.

In this paper, we extend the interior gradient estimate due to Korevaar-Simon
\cite{K} (see also \cite{K2}) for solutions of the prescribed mean curvature
equation in  Euclidean space to the general case of Killing graphs. Our
main application is the extension from $C^{2,\alpha}$ to $C^{0}$ boundary data 
of the result in \cite{DHL} on the solvability of the Dirichlet problem for 
the prescribed mean curvature equation, in particular, extending Theorem 16.11 of 
\cite{GT} to Killing graphs. In turn, this result has an interesting application 
in the case of ambient spaces of constant sectional curvature, in particular, 
for Plateau's problem for constant mean curvature hypersurfaces inside a solid 
torus in Euclidean space; see our Theorem \ref{pla} below. 

In addition to the above, we use our gradient estimate and a result in \cite{DHL}
on the Dirichlet problem on bounded domains for
$C^{2,\alpha}$ boundary data to show the solvability of the
asymptotic Dirichlet problem for radial graphs in the hyperbolic space
with prescribed mean curvature function. A similar result for constant mean
curvature was obtained in \cite{GS} using a different approach. 
We observe that an interior gradient estimate for Killing graphs was
recently used in \cite{LW} to show the existence of nonparametric solutions for a
capillary problem in warped products. Other applications
of the Korevaar-Simon method can be found in \cite{EH} and \cite{Ei}.
\vspace{1ex}

To give the precise statements of our results we first discuss the geometric
setting. \vspace{1,5ex}

Let $N^{n+1}$ denote an $(n+1)$-dimensional Riemannian manifold
carrying a non-singular complete Killing vector field $Y$ whose orthogonal
distribution is integrable. Fix a integral hypersurface $M^n$ of the
orthogonal distribution. Observe that $M^n$ is a totally geodesic
submanifold of $N^{n+1}$. Let $\Omega\subset M^n$ be a domain such that the flow
 $\Psi\colon\R\times M^n\to N^{n+1}$ generated
by $Y$ has complete orbits. Notice that $\gamma=1/\<Y,Y\>$ can be seen as a function
in $\bar\Omega$ since $Y\gamma=0$ by the Killing equation. Moreover, the solid
cylinder $\Psi(\R\times\bar\Omega)$ with the induced metric has a
warped product Riemannian structure $\bar\Omega\times_{\rho}\R$ where
$\rho=1/\sqrt{\gamma}$.

Given a function $u$ on $\bar\Omega$ the associated \emph{Killing graph} is
the hypersurface
$$
\mbox{Gr}(u)=\{\Psi(u(x),x): x\in\bar\Omega\}.
$$
It was shown in \cite{DHL} that $\mbox{Gr}(u)$ has mean curvature $H(x)$ if
and only if $u\in C^2(\Omega)$ satisfies
\be\label{pde}
\mathcal{Q}[u]=\text{div}\left(\frac{\nabla u}{w}\right) 
-\frac{\gamma}{w}\<\nabla u, \bar\nabla_YY\>=nH
\ee
where $w=\sqrt{\gamma+|\nabla u|^2}$ and $H$ is computed for the orientation
given by the Gauss map
\be\label{gauss}
N=\frac{1}{w}(\gamma Y-\Psi_{*}\nabla u)
\ee
where $\nabla$ and $\text{div}$ denote the gradient and divergence in $M^n$
and $\bar\nabla$ the Riemannian covariant derivative in $N^{n+1}$.
\vspace{1,5ex}

Given a point $o\in\Omega$, let $r>0$ be such that $r<i(o)$ where $i(o)$ is the
injectivity radius of $M^n$ at $o$. We denote by $B_{r}(o)$ the geodesic
ball contained in $\Omega$ centered at $o$ and radius~$r$. We now state
our interior gradient estimate for Killing graphs. A similar result but just 
for product ambient spaces $N^{n+1}=M^n\times\R$ was given in \cite{S}.

\begin{theorem}\label{korevaar}\po Let $u\in C^3(B_r(o))$ be a
negative solution of the mean curvature equation~(\ref{pde}) for $H\in
C^1(\Omega)$. Then, there exists a constant $L=L(u(o),r,\gamma,H)$ such that
$|\nabla u(o)|\leq L$.
\end{theorem}

The \emph{Killing cylinder} $K$ over $\Gamma=\partial\Omega$ is the hypersurface
in $N^{n+1}$ ruled by the flow lines of $Y$ through $\Gamma$ given by 
$$
K=\{q=\Psi(s,x):s\in\R\;\mbox{and}\;x\in\Gamma\}.
$$
In the sequel, the mean curvature of $K$ pointing inward is denoted
by $H_{\textrm{cyl}}$.\vspace{1,5ex}

The use of Theorem \ref{korevaar} gives the following extension of the main result in \cite{DHL}.

\begin{theorem}\label{dhl}\po 
Let $\Omega\subset M^n$ be a $C^{2,\alpha}$ bounded domain  satisfying $H_{cyl}\geq 0$ and
$$
\inf_{\Omega}\operatorname*{Ric}_M\geq-n\inf_{\Gamma}H_{cyl}^{2}.
$$
Let $H\in C^\beta(\Omega)$ for some $0<\beta<1$ be such that $\sup_\Omega|H|\leq\inf_\Gamma H_{cyl}$. 
Given $\varphi\in C^{0}(\Gamma)$ there exists
a unique function $u\in C^{2,\beta}(\Omega)\cap C^0(\bar{\Omega})$ whose
Killing graph has mean curvature $H$ and $u|_\Gamma=\varphi$.
\end{theorem}

 Similar extensions of Theorems $2$ and $3$ of \cite{DHL} (proved under weaker assumptions 
on the Ricci curvature of $M^n$ than in Theorem $1$ in that paper) also follow. 
Observe also that Theorem \ref{dhl} generalizes Theorem 1.4 of \cite{S}.
\vspace{1ex}

Although our setting is very general and the proof is
based on nonparametric methods, it turns out that our Theorem \ref{dhl} has an interesting 
consequence on Plateau's problem for constant mean curvature hypersurfaces in $\R^{n+1}$, 
the hyperbolic space $\Hy^{n+1}$ and the sphere $\Sp^{n+1}$. Let $\Q^{n+1}$ stand for any of
these spaces and let $\Q_+^n$ be a half totally geodesic hypersurface of $\Q^{n+1}$, hence 
$\Q_+^n$ is complete up to the boundary and $\partial\Q_+^n$ is a totally geodesic 
submanifold of  $\Q^{n+1}$ of codimension two. We denote by $T$ the  solid torus in
$\Q^{n+1},n\geq2$, obtained by rotating a $C^{2,\alpha}$ domain
$\Omega\subset \Q_+^n\setminus\partial\Q_+^n$ diffeomorphic to a ball 
around $\partial\Q_+^n$ and by $H_{\cal{T}}$ the mean curvature of its boundary
$\cal{T}$ with respect to the inner orientation.

\begin{theorem}\po\label{pla} Assume that $\Omega$ is $C^{2,\alpha}$ and that $H_{\cal{T}}\geq H$
for some constant $H\geq 0$. If $\Gamma\subset\cal{T}$ is a compact embedded topological 
hypersurface of $\cal{T}$ that intersects at a single point every circle through $\partial\Omega$ 
orthogonal to $\Q_+^n$, then there exists a unique compact embedded topological hypersurface
\mbox{$M^n\subset T$} with boundary $\Gamma$ homeomorphic to a closed ball 
that is~$C^{\infty}$ and has constant mean curvature $H$ on $M^n\setminus\Gamma$.
\end{theorem}

Observe that the mean curvature  with respect to the inner
orientation of a geodesic sphere in $\Q_+^n\setminus\partial\Q_+^n$ 
with radius $r$  tends to $+\infty$ as $r\rightarrow 0$. Thus, the mean curvature 
of the boundary of the solid torus obtained by rotating the geodesic ball in
$\Q_+^n\setminus\partial\Q_+^n$ has mean curvature larger than any given $H>0$ 
if the radius of the ball is small enough, hence  Theorem \ref{pla} applies to
such a torus.
\vspace{1ex}

In this paper,  the other application of  Theorem \ref{korevaar} is to prove the
solvability of the asymptotic Dirichlet problem for radial graphs in the
hyperbolic space with prescribed mean curvature function.
Such a result but for constant mean curvature was obtained in \cite{GS} by 
a different approach. Our solution to the problem is obtained as the limit of radial
graphs with prescribed mean curvature on a sequence of concentric geodesic
balls exhausting $\Hy^n$. The existence of these Killing graphs on
geodesic balls follows from the general result in \cite{DHL}. \vspace{1ex}

Let $\Hy^n$ be a complete totally geodesic hypersurface in
$\Hy^{n+1}$. A fixed constant speed geodesic line $\ell$ orthogonal to
$\Hy^n$ at a point $o\in\Hy^n$ determines uniquely a
one-parameter family of translation isometries in $\Hy^{n+1}$ that
preserve $\ell$ and whose Killing field $Y$ extends the velocity vector of
$\ell$. These isometries extend to the ideal boundary at infinity
$\partial_\infty\Hy^{n+1}$ of $\Hy^{n+1}$ as conformal
transformations belonging to the conformal structure of $\partial_\infty\Hy^{n+1}$.

Denote by $\Psi\colon\,\R\times\Hy^n\to\Hy^{n+1}$ the
flux generated by $Y$ and set
$\bar\Hy^n=\Hy^n \cup\partial_\infty\Hy^n$.
Given a function $u\in C^2(\Hy^n)\cap C^{0}(\bar\Hy^n)$
its \emph{radial graph} (Killing graph) is the hypersurface
$$
\mbox{Gr}(u) =\{\Psi(u(x),x): x\in\Hy^n\}
$$
with asymptotic boundary given by
$$
\partial_\infty\mbox{Gr}(u)=\overline{\mbox{Gr}(u)}\cap\partial_\infty\Hy^{n+1}.
$$
With some abuse of language, we say that $\mbox{Gr}(\phi)=\{\Psi(\phi(x),x):
x\in\partial_\infty\Hy^n\}$ is the radial graph of the function
$\phi= u|_{\partial_\infty\Hy^n}$. \medskip

\begin{theorem}\label{thm1}\po 
 Given $H\in C^{\beta}(\Hy^{n})$, $0<\beta<1$, with  $\sup_{\Hy^{n}}|H|<1$  
and $\phi\in C^0(\partial_\infty\Hy^{n})$ there exists a unique function $u\in
C^{2,\beta}(\Hy^{n})\cap C^{0}(\bar{\Hy}^{n})$ such that $\mbox{Gr}(u)$ has mean 
curvature $H$ and $\partial_\infty\mbox{Gr}(u)$ is the radial graph of $\phi$.
\end{theorem}

The above result is a natural extension to Killing graphs in the hyperbolic
space of the asymptotic Dirichlet problem for the mean curvature hypersurface
PDE,  which in the minimal case has been studied in a complete noncompact 
Riemannian manifold in \cite{EFR}, \cite{GR} and \cite{RT}.


\section{The gradient estimate}

For a system of coordinates $x=(x^{1},\ldots, x^{n})$ equation (\ref{pde})
reads as
\be\label{pde2}
a^{ij}u_{i;j}-R\<\nabla\gamma,\nabla u\>=nHw
\;\;\;\;\mbox{where}\;\;\;a^{ij}=\sigma^{ij}-\frac{u^iu^j}{w^2}.
\ee
Here the $\sigma_{ij}$'s are the coefficients of the metric in $M^n$ and for
simplicity we denote
$$
R=\frac{\gamma+ w^2}{2\gamma w^2}.
$$

\noindent\emph{Proof of Theorem \ref{korevaar}:}
We define a nonnegative function $\eta(x)=g(\phi(x))$ on $B_{r}(o)$ where
$$
g(t) = e^{C_{1}t}-1
$$
for some constant $C_{1}>0$ and $\phi$ is given by
$$
\phi(x) = \left(  1-\frac{d^2(x)}{r^2}+\frac{u(x)}{2u_0}\right)  ^{+}.
$$
Here $+$ means positive part, $d(x)$ is the geodesic distance to $o$ in
$M^n$ and $C=1/2u_{o}$ where $u_0=-u(o)$. Then $\eta$ vanishes outside of
$B_{r}(o)$. To be more precise, we should replace $1$ by $1-\epsilon$ in the
definition of $\phi$ so that $\eta$ is smooth with compact support and later
let $\epsilon$ tend to zero, but this is omitted for simplicity.

Let $p\in B_{r}(o)$ be an interior point where the function $h=\eta w$ has a
maximum. In the sequel computations are done at this point without further
notice. From $h_i=0$ we have
\be\label{one}
\eta_i w =-\eta w_i.
\ee
Moreover, the Hessian matrix of $h$ is negative semidefinite. Taking the trace
of the product of the Hessian matrix of $h$ with the positive definite matrix
$a^{ij}/w$ yields
$$
0\ge\frac{1}{w} a^{ij}h_{i;j}=\frac{a^{ij}}{w}(w\eta_{i;j} + 2\eta_iw_j
+\eta w_{i;j}).
$$
We obtain using (\ref{one}) that
\be\label{five}
a^{ij}\eta_{i;j} +\frac{\eta}{w^2}a^{ij}(ww_{i;j}-2w_iw_j)\leq 0.
\ee
 From (\ref{gauss}) we have
\be\label{en}
N^k= -\frac{u^k}{w}.
\ee
Thus,
\be\label{three}
w_i=\frac{\gamma_i}{2w}+\frac{u^ku_{k;i}}{w}=\frac
{\gamma_i}{2w}-N^ku_{k;i}.
\ee
In the sequel, we use (\ref{one}), (\ref{en}) and (\ref{three}) several times
without further reference.

We have,
\bea
w_{i;j}\!\!\!  &  =&\!\!\! \frac{\gamma_{i;j}}{2w}-\frac{\gamma_iw_j
}{2w^2}-N^k_{;j} u_{k;i} -N^ku_{k;ij}\\
\!\!\!  &  =&\!\!\! \frac{\gamma_{i;j}}{2w}+\frac{\sigma^{kl}}{w}u_{l;j}u_{k;i}
-\frac{w}{\eta^2}\eta_i\eta_j-N^k u_{k;ij}\\
\!\!\!  &  =&\!\!\! \frac{\gamma_{i;j}}{2w} +\frac{1}{w}\big(\sigma^{kl}-N^k
N^l\big)u_{l;j}u_{k;i}+\frac{1}{w}N^k N^l u_{l;j}u_{k;i} -\frac{1}
{w}w_iw_j-N^k u_{k;ij}\\
\!\!\!  &  =&\!\!\! \frac{\gamma_{i;j}}{2w} +\frac{a^{kl}}{w} u_{l;j}
u_{k;i}+\frac{\gamma_i\gamma_j}{4w^{3}} -\frac{1}{2w^2}(w_i\gamma
_j+w_j\gamma_i)-N^ku_{k;ij}.
\eea
Thus,
\be\label{main}
a^{ij}w_{i;j}=\frac{1}{2w}a^{ij}\gamma_{i;j} +\frac{1}{w}
a^{ij}a^{kl}u_{l;j}u_{k;i}+\frac{1}{4w^{3}}a^{ij}\gamma_i\gamma_j
+\frac{1}{w\eta}a^{ij}\eta_i\gamma_j -N^ka^{ij}u_{k;ij}.
\ee
We compute the last term of (\ref{main}). The Ricci identities for the Hessian
of $u$ yield
$$
u_{k;ij} = u_{i;kj}=u_{i;jk}+R^l_{kji}u_l.
$$
Hence,
$$
N^k a^{ij}u_{k;ij} =N^ka^{ij} u_{i;jk}- \frac{1}{w}a^{ij}R_{kji}^l u^k
u_l =N^k (a^{ij}u_{i;j})_{;k}-N^k a^{ij}_{;k}u_{i;j} -\frac{1}{w}
a^{ij}R_{kji}^l u^k u_l.
$$
It follows that
\be\label{equation}
N^k a^{ij}u_{k;ij}= nN^k (wH)_k + N^k(R\<\nabla
\gamma,\nabla u\>)_k - N^ka^{ij}_{;k}u_{i;j}-\frac{1}{w}a^{ij}R_{kji}^l
u^k u_l.
\ee
We compute the first three terms of (\ref{equation}). For the first term, we
have
\be\label{first}
(wH)_k = \frac{w}{\eta}(\eta H_k-H\eta_k).
\ee
Since
$$
\Big(\frac{\gamma+w^2}{w^2}\Big)_k =\frac{\gamma_k}{w^2}
-\frac{\gamma}{w^{4}}(\gamma_k+2u^lu_{l;k}) =\frac{1}{w^2}
\Big(\gamma_k+\frac{2\gamma\eta_k}{\eta}\Big)
$$
and
$$
\Big(\frac{1}{2\gamma}\<\nabla\gamma,\nabla u\>\Big)_k =\Big(\frac
{\gamma_l}{2\gamma}\Big)_{;k}u^l +\frac{\gamma^l}{2\gamma}u_{l;k}
=\frac{1}{2\gamma}\left[  \Big(\frac{\gamma_l\gamma_k}{\gamma}
-\gamma_{k;l}\Big)wN^l +\gamma^lu_{l;k}\right]  ,
$$
we obtain for the second term of (\ref{equation}) that
\be\label{second}
(R\<\nabla\gamma,\nabla u\>)_k=R\Big[\Big(\frac{\gamma
_l\gamma_k}{\gamma}-\gamma_{k;l}\Big)wN^l +\gamma^lu_{l;k}
\Big]+\frac{1}{w^2}\Big(\frac{\gamma_k}{2\gamma} 
+\frac{\eta_k}{\eta}\Big)\<\nabla\gamma,\nabla u\>.
\ee
We have,
\bea
a^{ij}_{;k} \!\!\!  &  =&\!\!\! -\frac{1}{w^2}(u^i_{;k}u^j+u^i
u^j_{;k})+\frac{1}{w^{4}}(\gamma_k -2wN^lu_{l;k})u^i u^j\\
\!\!\!  &  =&\!\!\!\frac{1}{w}(u^i_{;k}-N^i N^l u_{l;k})N^j 
+ \frac{1}{w}(u^j_{;k}-N^j N^l u_{l;k})N^i+\frac{1}{w^2}\gamma_k
N^iN^j\\
\!\!\!  &  =&\!\!\!\frac{1}{w}a^{il}u_{l;k} N^j + \frac{1}{w}a^{jl}
u_{l;k}N^i+\frac{1}{w^2}\gamma_kN^iN^j.
\eea
It follows that the third term of (\ref{equation}) is given by
\be\label{third}
N^k a^{ij}_{;k}u_{i;j} =\frac{2}{w}a^{ij}\Big(\frac{w\eta_i
}{\eta}+\frac{\gamma_i}{2w}\Big) \Big(\frac{w\eta_j}{\eta}+\frac
{\gamma_j}{2w}\Big)+\frac{1}{w^2}\gamma_kN^kN^i 
\Big(\frac{w\eta_i}{\eta}+\frac{\gamma_i}{2w}\Big).
\ee
Replacing (\ref{first}), (\ref{second}) and (\ref{third}) into (\ref{equation}
) yields
\bea
N^ka^{ij}u_{k;ij}\!\!\!&=&\!\!\!n\frac{w}{\eta}N^k(\eta H_k-H\eta
_k) -\frac{2}{w}a^{ij}\Big(\frac{w\eta_i}{\eta}+\frac{\gamma_i}
{2w}\Big)\Big(\frac{w\eta_j}{\eta}+\frac{\gamma_j}{2w}\Big)\\
\!\!\!&-&\!\!\!\frac{1}{w^2}\gamma_kN^kN^i \Big(\frac{w\eta_i
}{\eta}+\frac{\gamma_i}{2w}\Big) +\frac{1}{w^2}N^k\Big(\frac{\gamma_k
}{2\gamma} +\frac{\eta_k}{\eta}\Big)\<\nabla\gamma,\nabla u\>\\
\!\!\!&+&\!\!\!R \Big[\Big(\frac{\gamma}{2w}\sigma^{kl}+wN^k
N^l\Big) \frac{\gamma_k \gamma_l}{\gamma}-wN^k N^l\gamma_{k;l}
+\frac{w}{\eta}\gamma_l\eta^l\Big]-\frac{1}{w}a^{ij}R_{kji}^l u^k
u_l.
\eea
It now follows from (\ref{main}) that
\bea
&&a^{ij}w_{i;j}-\frac{2}{w}a^{ij}w_i w_j = \frac{3}{4w^{3}}a^{ij}
\gamma_i\gamma_j +\frac{1}{w}a^{ij}a^{kl}u_{l;j}u_{k;i} +\frac{3}{w\eta
}a^{ij}\gamma_i\eta_j+\frac{1}{2w}a^{ij}\gamma_{i;j}\\
&& +\frac{1}{w}a^{ij}R_{kji}^l u^ku_l -nN^k\frac{w}{\eta}(\eta
H_k-H\eta_k) +\frac{1}{w^2}\gamma_kN^kN^i \Big(\frac{w\eta_i
}{\eta}+\frac{\gamma_i}{2w}\Big)\\
&&  -\frac{1}{w^2}N^k\Big(\frac{\gamma_k}{2\gamma} +\frac{\eta_k}{\eta
}\Big)\<\nabla\gamma,\nabla u\> -R\Big[\Big(\frac{\gamma}{2w}\sigma^{kl}
+wN^kN^l\Big)\frac{\gamma_k \gamma_l}{\gamma}-wN^k N^l\gamma_{k;l}
+\frac{w}{\eta}\gamma_l\eta^l\Big].
\eea
After multiplying both sides by $\eta/w$ and discarding the non-negative
terms, we obtain
\bea
&&  \frac{\eta}{w}\Big(a^{ij}w_{i;j}-\frac{2}{w}a^{ij}w_i w_j
\Big) \geq\Big[-nN^kH_k-\frac{\gamma_k}{2\gamma w^{3}}N^k
\<\nabla\gamma,\nabla u\>+\frac{1}{2w^2}a^{ij}\gamma_{i;j}\\
&&  -R\Big(\Big(\frac{\gamma}{2w^2}\sigma^{kl} +N^k N^l
\Big)\frac{\gamma_k \gamma_l}{\gamma} -N^kN^l\gamma_{k;l}
\Big) +\frac{1}{w^2} a^{ij}R_{kji}^l u^ku_l\Big]\eta\\
&&+\Big[\Big(nH +\frac{1}{w^2}N^k\gamma_k -\frac{1}{w^{3}
}\<\nabla\gamma,\nabla u\> \Big)N^i+\Big(\frac{3}{w^2}a^{ij} -R\sigma
^{ij}\Big)\gamma_j\Big]\eta_i.
\eea
It is easy to check that there is a positive constant $M=M(\gamma,H)$ such
that
\be\label{six}
\frac{1}{w^2}\eta a^{ij}(ww_{i;j}-2w_i w_j)\ge-M\eta
-A^i\eta_i
\ee
where $-A^i$ denotes the coefficient of $\eta_i$. It follows from
(\ref{five}) and (\ref{six}) that
\be\label{des}
a^{ij}\eta_{i;j}-M\eta-A^i\eta_i\leq0.
\ee
On the other hand,
$$
\eta_i = g^{\prime-2}(d^2)_i +C u_i)
$$
and
$$
\eta_{i;j}= g^{\prime-2}(d^2)_{i;j} +C u_{i;j})+g^{\prime\prime-2}
(d^2)_i +C u_i)(-r^{-2}(d^2)_j +C u_j).
$$
Hence,
\bea
a^{ij}(r^{-2}(d^2)_i \!\!\!&-&\!\!\! C u_i)(r^{-2}(d^2)_j-Cu_j)\\
\!\!\!&=&\!\!\!\frac{C^2\gamma}{w^2}|\nabla u|^2 -\frac{2C\gamma
}{r^2w^2}\<\nabla u,\nabla d^2\> +\frac{1}{r^{4}} \Big(|\nabla
d^2|^2-\frac{1}{w^2}\<\nabla u, \nabla d^2\>^2 \Big)\\
\!\!\!&\geq&\!\!\!\frac{C^2\gamma}{w^2}\Big(|\nabla u|^2 
-\frac{2}{Cr^2}\<\nabla u,\nabla d^2\>\Big)
\eea
and
$$
a^{ij}(-r^{-2}(d^2)_{i;j} +C u_{i;j}) = -r^{-2}a^{ij}(d^2)_{i;j}
+C(nHw+R\<\nabla\gamma,\nabla u\>)
$$
where
$$
a^{ij}(d^2)_{i;j}=\Delta d^2 -\frac{1}{w^2}\<\nabla_{\nabla u}\nabla
d^2,\nabla u\>.
$$
It follows from (\ref{des}) that
\bea
\frac{C^2\gamma}{w^2}\Big[|\nabla u|^2 
\!-\!\frac{2}{Cr^2}\<\nabla u,\nabla d^2\>\Big]g^{\prime\prime}
\!\!\!\!\!  &+&\!\!\!\!\!\Big[\frac{1}{r^2w^2}\<\nabla_{\nabla u}\nabla d^2,\nabla u\>
-\frac{\Delta d^2}{r^2}\!+\! C(nHw+R\<\nabla\gamma,\nabla u\>)\Big]g^{\prime}\\
\!\!\!&\leq&\!\!\! Mg+A^i(-r^{-2}(d^2)_i+Cu_i)g^{\prime}.
\eea
Since
$$
-A^i(d^2)_i=\Big(nH +\frac{1}{w^2}N^k\gamma_k 
-\frac{1}{w^3}\<\nabla\gamma,\nabla u\>\Big)N^i (d^2)_i +\Big(\frac{3}{w^2}a^{ij}
-R\sigma^{ij}\Big)\gamma_j(d^2)_i
$$
and
$$
A^i u_i = \Big(nHw +\frac{1}{w}N^k\gamma_k -\frac{1}{w^2}
\<\nabla\gamma,\nabla u\>\Big)
\frac{|\nabla u|^2}{w^2}+\frac{3}{w}a^{ij}\gamma_jN_i +R\<\nabla
\gamma,\nabla u\>,
$$
we conclude that
$$
\frac{C^2\gamma}{\gamma+|\nabla u|^2}\Big(|\nabla u|^2 -\frac{2}{Cr^2
}\<\nabla u,\nabla d^2\>\Big)g^{\prime\prime}+Pg^{\prime}-Mg \leq0
$$
where the coefficient
\bea
\!\!\!&&\!\!\!P=\frac{n}{w}CH\gamma-\frac{1}{r^2}\Big(\Delta d^2
-\frac{1}{w^2}\<\nabla_{\nabla u}\nabla d^2,\nabla u\> \Big)-C\Big(wN^k
\gamma_k-\<\nabla\gamma,\nabla u\> \Big)\frac{|\nabla u|^2}{w^{4}}\\
\!\!&-&\!\!\!\!\!\frac{3}{w}Ca^{ij}\gamma_jN_i-\frac{1}{r^2w^{3}
}\Big(nHw^{3}+wN^k\gamma_k -\<\nabla\gamma,\nabla u\>\Big)N^i
(d^2)_i -\frac{1}{r^2}\Big(\frac{3}{w}a^{ij} -R\sigma^{ij}
\Big)\gamma_i (d^2)_j
\eea
of $g^{\prime}$ is bounded above by a positive constant $C_0=C_0
(u(o),r,\gamma, H)$.

Suppose that
$$
|\nabla u|\ge\frac{8}{Cr}=\frac{16u_0}{r}.
$$
Then, we have
$$
|\nabla u|\ge\frac{4|\nabla d^2|}{Cr^2}
$$
and
$$
|\nabla u|^2-\frac{2}{Cr^2}\<\nabla u,\nabla d^2\>\ge|\nabla u|^2
-\frac{2}{Cr^2}|\nabla u||\nabla d^2|\ge\frac{1}{2}|\nabla u|^2.
$$
Thus, there exists a constant $D=D(u(0),r,\gamma)$ such that
$$
\frac{C^2\gamma}{\gamma+|\nabla u|^2}\Big(|\nabla u|^2 -\frac{2}{Cr^2
}\<\nabla u,\nabla d^2\>\Big)
\geq\frac{C^2\gamma|\nabla u|^2}{2(\gamma+|\nabla u|^2)}\geq D>0.
$$
For instance, setting $\gamma_0=\inf_{B_{r}(o)}\gamma$ we may take
$D=32\gamma_0/(r^2\gamma_0+256u^2(0)). $ Since $g(t) = e^{C_{1}t}-1$,
we obtain a contradiction taking $C_{1}=C_{1}(u(o),r,\gamma, H)$ sufficiently
large, that is, such that
$$
Dg^{\prime\prime}(p)+C_0g^{\prime}(p)-Mg(p)>0.
$$
Hence,
$$
w(p)\leq C_2=\gamma_{1}+\frac{16u_0}{r}
$$
where $\gamma_{1}=\sup_{B_{r}(o)}\gamma$. Since $h(o)\leq h(p)$, we obtain
$\eta(o)w(o)\leq\eta(p)w(p)$.
Therefore,
$$
(e^{C_{1}/2}-1)w(o)\le C_2e^{C_{1}}.
$$
This gives the desired estimate and concludes the
proof.\qed

\begin{remark}\po \emph{Observe that the constant $L$ given by Theorem
\ref{korevaar} also depends on the geometry of $M^n$ along $B_{r}(o)$
including the Ricci curvature. }
\end{remark}

\section{The applications}

In this section, we provide the proofs of the applications of Theorem \ref{korevaar}
that has been stated in the introduction.
\vspace{1,5ex}

\noindent\emph{Proof of Theorem \ref{dhl}:} Let $\varphi_k^\pm\in
C^{2,\alpha}(\Gamma)$ be monotonic sequences converging from above and below
to $\varphi$ in the $C^0$ norm. Let $H_k\in C^{\infty}(\Omega)$ 
be a sequence converging to $H$ such that 
$H_0=\sup_k\sup_{\Omega}|H_k|\leq\inf_{\Gamma}H_{cyl}$. Set $H^+=H_0$ and
$H^-=-H_0.$ By Theorem 1 of \cite{DHL} there exist solutions 
$u_k^\pm\in C^{2,\alpha}(\bar{\Omega})$ of $\mathcal{Q}_{H_k}[u]=nH_k$ and
$v_k^\pm\in C^{2,\alpha}(\bar{\Omega})$ of $\mathcal{Q}_{H^\pm}[v]=nH^\pm$ 
such that $u_k^\pm|_{\Gamma}=\varphi_k^\pm$ and $v_k^\pm|_{\Gamma}=\varphi_k^\pm$. 
>From Theorem 6.17 of \cite{GT} it follows that
$u_k^\pm,v_k^\pm\in C^{\infty}(\Omega)$. The sequences $u_k^\pm$,
$v_k^\pm$ have uniformly bounded $C^0$ norm since, by the comparison
principle,
\be\label{comp1}
v_1^-\leq...\leq v_k^-\leq v_{k+1}^-\leq...\leq v_{k+1}^+\leq v_k^+
\leq...\leq v_1^+,\;\;\;k=1,2,...
\ee
and 
\be\label{comp2}
v_k^-\leq u_k^\pm\leq v_k^+,\;\;\;k=1,2,...
\ee
It follows from Theorem \ref{korevaar} and linear elliptic PDE
theory that the sequences $u_k^\pm$, $v_k^\pm$ have equicontinous
$C^{2,\beta}$ norm on compact subsets of $\Omega$, the uniform $C^{2,\beta}$
norm depending only on the distance of the compact subset to $\Gamma$
(Corollary 6.3 of \cite{GT}). Considering an exhaustion of $\Omega$ by
relatively compact subdomains and using the diagonal method we obtain that
$u_k^\pm$ and $v_k^{\pm}$ contain subsequences converging uniformly on
compact subsets of $\Omega$ on the $C^2$ norm to $u,v^{\pm}\in C^2(\Omega)$ 
satisfying $\mathcal{Q}_{H}[u]=nH$ and $\mathcal{Q}_{H^\pm}[v^\pm]=nH^\pm$. 
By (\ref{comp2}) we have that $v^{-}\leq u\leq v^+$. 
Since $v_k^\pm|_{\partial\Omega}=\varphi_k^\pm$ converges to $\varphi$, 
it follows from (\ref{comp1}) that $v^\pm$ extends
continuously to $\bar{\Omega}$ and $v^\pm|_{\partial\Omega}=\varphi$. Therefore,
$u$ extends continuously to $\bar{\Omega}$ and $u|_{\partial\Omega}=\varphi$.
Since $H\in C^\beta(\Omega)$ and $\mathcal{Q}_{H}[u]=nH$ 
then $u\in C^{2,\beta}(\Omega)\cap C^0(\bar\Omega)$ from 
Theorem 6.17 of \cite{GT}.
\qed\vspace{1,5ex}

\noindent\emph{Proof of Theorem \ref{pla}: }The circles orthogonal to
$\Q_+^n$ are the orbits of a Killing field  orthogonal to
$\Q_+^n$. By assumption $\Gamma$ is the graph of a continuous
function $\phi\in C^{\infty}\left(  \partial\Omega\right)  $. It is clear
that the other conditions of Theorem \ref{dhl} are satisfied when applied 
to the universal cover of $T$, and thus the proof
follows.\qed\vspace{1,5ex}

\noindent\emph{Proof of Theorem \ref{thm1}:}   The proof will be done in two steps. 
First, we show the existence of a uniform height estimate for the solutions of 
the Dirichlet problem on a sequence of concentric geodesic balls exhausting $\Hy^n$. 
Then, we prove that a subsequence of these solutions converges to a
solution of our problem.  We may assume that $H\in C^\infty(\Hy^n)$ since  
the case $H\in C^\beta(\Hy^n)$ can be obtained by using approximating arguments as in the
proof of Theorem \ref{dhl}.\vspace{1ex}

We represent the hyperbolic space $\Hy^{n+1}$ of constant sectional
curvature $-1$ as the warped product manifold
$$
\Hy^{n+1}=\Hy^n\times_{\cosh\rho}\R
$$
where $\rho$ is the geodesic distance in $\Hy^n$ to the point
$o\in\Hy^n$ determined by the geodesic $\ell$. In terms of the
notation fixed earlier $Y=\partial/\partial s$ where $s$ parametrizes the
factor $\R$. In (\ref{pde}) we thus have
\be\label{gamma}
\gamma=1/\cosh^2\rho.
\ee

Let $\Omega_k$ denote the geodesic disc in $\Hy^n$ of radius
$\rho_k$ centered at $o$ with boundary $\Gamma_k$ where
$$
\rho_k=\text{arctanh}(1-1/k),\;\;k=2,3,\ldots
$$
Let $F\in C^{2,\alpha}(\Hy^n)\cap C^{0}(\bar\Hy^n)$ be such
that $F|_{\partial_\infty\Hy^n}=\phi$ and set $M_0=|F|_0$. We
claim that the unique solution $u_k\in C^2(\bar\Omega_k)$ for any
$k\ge2$ obtained in \cite{DHL} of the auxiliary Dirichlet problem
\bea\left\{\begin{array}
[c]{cc}
\!\!\mathcal{Q}[u_k]=nH|_{\Omega_k}\vspace{1ex} & \\
\!\!\!\!\!u_k|_{\Gamma_k}=F|_{\Gamma_k}, &
\end{array}\right.
\eea
has a uniform height estimate, that is, independent of $k$.

The geodesic distance $d$ from $\Gamma_k$ measured along normal
geodesics pointing towards $o$ is $d=\rho_k-\rho$.
The geodesic curvature of the flow lines of $Y$ is given by
\be\label{kappa}
\kappa=\<\nabla d, \bar\nabla_{\sqrt{\gamma}Y} \sqrt{\gamma}Y\>
=-\gamma\<Y,\bar\nabla_{\nabla d}Y\> =\frac{\gamma^{\prime}}{2\gamma}
=\tanh\rho
\ee
where $^{\prime}$ denotes derivative with respect to $d$ and we used
(\ref{gamma}). Since the mean curvature of the geodesic ball $\Gamma_k$ 
is $\coth\rho_k$, then the constant mean curvature $H_k$ of the Killing 
cylinder  over $\Gamma_k$ is
\be\label{hk}
H_k=\frac{1}{n}((n-1)\coth\rho_k+\tanh\rho_k)>1
\ee
with respect to the pointing inward orientation.

Consider on $\Omega_k$ the function $v_k(x) = M_0 + h(d(x))$
where $h$ has to be chosen. Then,
\bea
\mathcal{Q}[v_k] \!\!\!&=&\!\!\! \text{div}\Bigg(\frac{h^{\prime}\nabla
d}{\sqrt{\gamma+h^{\prime2}}}\Bigg)- \frac{\gamma h^{\prime}}{\sqrt
{\gamma+h^{\prime2}}} \<\nabla d,\bar\nabla_YY\>\\
\!\!\!&=&\!\!\! \frac{h^{\prime}}{\sqrt{\gamma+h^{\prime2}}}\big(\Delta d -
\kappa\big) + \Bigg(\frac{h^{\prime}}{\sqrt{\gamma+h^{\prime2}}}
\Bigg)^{\prime}.
\eea
Since $-\frac{1}{n-1}\Delta d$ is the mean curvature of the geodesic sphere in
$\Hy^n$ of radius $\rho=\rho_k-d$, we have that
\be\label{Hk}
\Delta d-\kappa=-nH_k.
\ee
Thus,
$$
\mathcal{Q}[v_k] = -\frac{nh^{\prime}H_k}{\sqrt{\gamma+h^{\prime2}}} +
\frac{\gamma(h^{\prime\prime}-\kappa h^{\prime})}{(\gamma+h^{\prime2}
)^{\frac{3}{2}}}.
$$

We choose
$$
h(d) = C(\arcsin(\tanh\rho_k)- \arcsin(\tanh\rho))
$$
where $C>0$ is a constant to be chosen and $\arcsin:[-1,1]\to[-\frac{\pi}
{2},\frac{\pi}{2}]$. Thus
$$
h^{\prime}=C/\cosh\rho\;\;\;\mbox{and}\;\;\;h^{\prime\prime}=C\tanh\rho
/\cosh\rho.
$$
It follows from (\ref{kappa}) that $h^{\prime\prime}-\kappa h^{\prime}=0$ and
using (\ref{gamma}) that
$$
\frac{h^{\prime}}{\sqrt{\gamma+h^{\prime2}}} =\frac{C}{\sqrt{1+C^2}}.
$$
Therefore,
$$
\mathcal{Q}[v_k] = -\frac{C}{\sqrt{1+C^2}}\, n H_k.
$$
Hence, in order to obtain that
$Q[v_k]< Q[u_k]=nH$
it suffices to choose $C$ such that
\be\label{above}
\frac{C}{\sqrt{1+C^2}}H_k>|H|.
\ee
By assumption $H_0=\sup|H|<1$. Using (\ref{hk}) it follows that
(\ref{above}) holds if we choose $C>0$ such that
$$
C > \sqrt{\frac{H_0^2}{1-H_0^2}}\ge\sqrt{\frac{H_0^2}{H_k^2-H_0^2}}
\ge\sqrt{\frac{H^2}{H_k^2-H^2}}.
$$

We have for $x\in\Gamma_k$ that $v_k(x)\ge F|_{\Gamma_k}(x)$.
We conclude that $v_k$ is an upper barrier to $u_k$, that is,
$u_k(x) \le v_k(x)$, $x\in\Omega_k$.
This implies that the functions $u_k$ have a uniform height estimate since
$v_k(x)\leq M=M_0 + C\pi$,
and proves the claim.\vspace{1ex}

It follows from the above height estimate and the gradient estimate in Theorem
\ref{korevaar} that for given $m\in\mathbb{N}$ the sequence 
$\{u_k|_{\bar{\Omega}_m}\}_{k>2m}$ has equibounded $C^{1}$ norm in 
$\bar{\Omega}_m$. To be able to apply Theorem \ref{korevaar} we have to translate the
solutions by $-M$. By elliptic theory this sequence has equibounded
$C^{2,\alpha}$ norm in $\Omega_m$. Then, there is a subsequence
$\{u_{k_j^{m}}\}_j$ of $\{u_k\}_{k>2m} $ converging on $\Omega_m$ in
the $C^2$ norm to a function $v_m\in C^2(\Omega_m)$. Denoting
$Q_{H}[u]=\mathcal{Q}[u]-nH$, we have that $Q_{H}[v_m]=0$. It follows from
PDE regularity (cf.\ \cite{GT}) that $v_m\in C^{2,\alpha}(\Omega_m)$.

Consider a subsequence $\{u_{k_j^{m+1}}\}_j$ of $\{u_{k_j^{m}}\}_j$
which converges in $\Omega_{m+1}$ to a function 
$v_{m+1}\in C^{2,\alpha}(\Omega_{m+1})$ satisfying $Q_{H}[v_{m+1}]=0$. 
After iterating this process,
we obtain a subsequence $\{u_{k_{1+i}^{m+i}}\}_i$ of $\{u_k\}_k $ that
converges uniformly on compact subsets of $\Hy^n$ in the $C^2$ norm
to a function $u\in C^2(\Hy^n)$ satisfying $Q_{H}[u]=0$. By
simplicity, we also denote this subsequence by $\{u_k\}_{k\in\mathbb{N}}$.
We claim that $u$ extends continuously to $\partial_\infty\Hy^n $
and that $u|_{\partial_\infty\Hy^n}=\phi$.

We first prove that $\partial_\infty\mbox{Gr}(u)\subset\mbox{Gr}(\phi)$ by
showing that if $p\in\partial_\infty\Hy^{n+1}\setminus
\mbox{Gr}(\phi)$ then $p\notin\partial_\infty\mbox{Gr}(u)$. Consider
$p\in\partial_\infty\Hy^{n+1}\setminus\mbox{Gr}(\phi)$. Since
$\mbox{Gr}(\phi)$ is compact and $p\notin\mbox{Gr}(\phi)$ there exists an
equidistant hypersurface $E$ of $\Hy^{n+1}$ such that $\partial_\infty E$ 
separates $p$ from $\mbox{Gr}(\phi)$, that is, $p$ and
$\mbox{Gr}(\phi)$ are in distinct open connected components of 
$\partial_\infty\Hy^{n+1}\setminus\partial_\infty E$. Moreover, since
$H_0=\sup_{\Hy^n}|H|<1$, we may assume that $E$ has constant mean
curvature $H_0$ with respect to the unit normal vector field pointing to the
connected component $U$ of $\Hy^{n+1}\setminus E$ whose asymptotic
boundary contains $\mbox{Gr}(\phi)$. Set $G_k=\mbox{Gr}(u_k)$. By the
convergence of $u_k|_{\Gamma_k}$ to $\phi$ there is $k_0$ such
that $\partial G_k\subset U$ and then $\partial G_k\cap E=\emptyset$ for
all $k\geq k_0$. By the tangency principle $G_k\cap E=\emptyset$, i.e.,
$G_k\subset U$ for all $k\geq k_0$. It follows that $p\notin
\partial_\infty\mbox{Gr}(u)$.

Now consider a sequence $x_k\in\Hy^n$ converging to $x\in
\partial_\infty\Hy^n$. By the compactness of $\bar{\Hy
}^{n+1} $ there is a subsequence $\Psi(u(x_{k_j}),x_{k_j})$ of
$\Psi(u(x_k),x_k)$ converging to $z\in\bar{\Hy}^{n+1}$. Since
$x_k$ diverges and $\Psi(u(x_{k_j}),x_{k_j})\in\mbox{Gr}(u)$ it follows
that $z\in\partial_\infty\mbox{Gr}(u)$. From what we proved above
$z\in\mbox{Gr}(\phi)$ and hence $z=\Psi(\phi(x_0),x_0)$ for some $x_0
\in\partial_\infty\Hy^n$. Since $u$ is globally bounded, the
sequence $\{u(x_{k_j})\}_j\subset\R$ is bounded and thus contains
a subsequence $\{u(x_{k_{j_i}})\}_i$ converging to some $s_0\in\R$. 
Being the extension of $\Psi_{s_0}$ to 
$\bar{\Hy}^{n+1}$ continuous, we obtain that
$$
z=\lim_i\Psi(u(x_{k_{j_i}}),x_{k_{j_i}})=\Psi(s_0,x)
$$
and therefore $\Psi(\phi(x_0),x_0)=\Psi(s_0,x)$. 
Since $\Psi\colon\R\times\partial_\infty\Hy^n\to\partial_\infty
\Hy^{n+1}$ is injective, it follows that $x=x_0$ and $s_0=\phi(x_0)$. 
Because the limits are the same for any convergent subsequence
considered, it follows that $u(x_k)\to\phi(x)$ as $k\to\infty$, and this
proves the claim.

Finally, uniqueness follows from the maximum principle applied to the
difference of two solutions.\qed

{\renewcommand{\baselinestretch}{1} \hspace*{-20ex}\begin{tabbing}
\indent \= Marcos Dajczer\\
\> IMPA \\
\> Estrada Dona Castorina, 110\\
\> 22460-320 -- Rio de Janeiro -- Brazil\\
\> marcos@impa.br\\
\end{tabbing}}

\vspace*{-5ex}

{\renewcommand{\baselinestretch}{1} \hspace*{-20ex}\begin{tabbing}
\indent \= Jorge Herbert de Lira\\
\> Departamento de Matematica, \\
\> Universidade Federal do Ceara,\\
\> Bloco 914 -- Campus do Pici\\
\> 60455-760 -- Fortaleza -- Ceara -- Brazil\\
\> jorge.lira@mat.ufc.br
\end{tabbing}}

\vspace*{-2ex}

{\renewcommand{\baselinestretch}{1} \hspace*{-20ex}\begin{tabbing}
\indent \= Jaime Ripoll \\
\> Departamento de Matematica \\
\> Univ. Federal do Rio Grande do Sul \\
\> Av. Bento Gon\c calves 9500\\
\> 91501-970 -- Porto Alegre -- RS -- Brazil\\
\> jaime.ripoll@ufrgs.br
\end{tabbing}}


\begin{thebibliography}{99}                                                                                               


\bibitem{GS} B. Guan and J. Spruck. \emph{Hypersurfaces of constant mean
curvature in hyperbolic space with prescribed asymptotic boundary at
infinity.} Amer. J. Math. \textbf{122} (2000), 1039--1060.

\bibitem{DHL} M. Dajczer, P. Hinojosa and J. H. de Lira. \emph{Killing graphs
with prescribed mean curvature.} Calc. Var. Partial Differential equations
\textbf{33} (2008), 231--248.

\bibitem{EH} K. Ecker and G. Huisken. \emph{Interior estimates for
hypersurfaces moving by mean curvature.} Invent. math. \textbf{105} (1991), 547--569.

\bibitem {Ei} M. Eichmair. \emph{The Plateau problem for marginally outer
trapped surfaces.} J. Diff. Geometry \textbf{83} (2009), 551--583.

\bibitem {EFR} N. do Esp\'{\i}rito-Santo, S. Fornari and J. Ripoll. \emph{The
Dirichlet problem for the minimal hypersurface equation in $M\times\R
$ with prescribed asymptotic boundary.} J. Math. Pures Appl. \textbf{93}
(2010), 204--221.

\bibitem{GR} J. G\'{a}lvez and H. Rosenberg. \emph{Minimal surfaces and
harmonic diffeomorphisms from the complex plane onto certain Hadamard
surfaces.} Amer. J. of Math. \textbf{132} (2010), 1249--1273.

\bibitem{GT} D. Gilbarg and N. Trudinger. \emph{Elliptic partial differential
equations of second order}. Springer Verlag, Berlin-Heidelberg, 2001.

\bibitem{K} N. Korevaar. \emph{An easy proof of the interior gradient bound
for solutions to the prescribed mean curvature problem.} Proc. Symp. Pure
Math. \textbf{45} (1986), 81--89.

\bibitem{K2} N. Korevaar. \emph{A priori interior gradient bounds for
solutions to elliptic Weingarten equations.} Ann. Inst. H.
Poincar\'e \ Anal. Non Lin\'eaire \textbf{4} (1987), 405--421.

\bibitem{LW} J. H. de Lira and G. Wanderley.  \emph{Existence of nonparametric 
solutions for a capillary problem in warped products.} http://arxiv.org/abs/1307.2871

\bibitem{RT} J. Ripoll and M. Telichevesky, \emph{On the asymptotic Dirichlet
problem for some divergence form quasi-linear elliptic PDE's.} To appear on
the Transactions of the AMS.

\bibitem{Se} J. Serrin, \emph{Gradient estimates for solutions of nonlinear
elliptic and parabolic equations.} Contributions to nonlinear functional
analysis (Proc. Sympos., Math. Res. Center, Univ. Wisconsin, Madison, Wis.,
1971), 565--601. Academic Press, New York, 1971.

\bibitem{S} J. Spruck. \emph{Interior gradient estimates and existence
theorems for constant mean curvature graphs in $M\times\R$.} Pure and
Applied Mathematics Quarterly \textbf{3} (2007), 785--800.
\end{thebibliography}
\end{document}